\documentclass[a4paper,10pt]{article}
\usepackage{natbib}
\usepackage{amsmath,amsfonts,amssymb,latexsym,wasysym,mathrsfs,verbatim}
\usepackage{amsthm}
\usepackage{pstricks,a4wide,setspace,multirow,graphicx} 
\usepackage[scriptsize]{subfigure}
\usepackage[scriptsize,bf]{caption}
\usepackage[all,2cell]{xy} \UseAllTwocells \SilentMatrices
\usepackage[english]{babel}

\usepackage{hyperref}

\theoremstyle{plain}

\newtheorem{proposition}{Proposition}[section]

\newtheorem{corollary}{Corollary}[section]
\theoremstyle{definition}
\newtheorem{definition}{Definition}[section]

\theoremstyle{remark}
\newtheorem{remark}{Remark}[section]

\theoremstyle{example}
\newtheorem*{example}{Example}

\author{Francesca Collet\\
{\small Dipartimento di Matematica} \\
{\small Alma Mater Studiorum Universit\`a di Bologna}\\
{\small Piazza di Porta San Donato 5; 40126 Bologna, Italy} \\
{\small e-mail: francesca.collet@unibo.it}\\[0.5cm]
Fabrizio Leisen\\
{\small School of Mathematics, Statistics and Actuarial Sciences} \\
{\small Cornwallis building, University of Kent; CT2 7NF Cantebury, United Kingdom}\\
{\small e-mail: fabrizio.leisen@gmail.com}\\[0.5cm]
Fabio Spizzichino\\
{\small Dipartimento di Matematica} \\
{\small Sapienza Universit\`a di Roma}\\
{\small Piazzale Aldo Moro 5; 00185 Roma, Italy } \\
{\small e-mail: fabio.spizzichino@uniroma1.it}}

\title{Merging exchangeable occupancy models: $\mathcal{M}^{(a)}$- models and relation with the maximum entropy principle}
\date{}

\begin{document}

\maketitle

\begin{abstract}
\noindent In this paper a new transformation of occupancy models, called {\it merging}, is introduced. In particular, it will be studied the effect of merging on a class of occupancy models that was recently introduced in \cite{CoLeSpSu13}. These results have an interesting interpretation in the so-called {\it entropy maximization inference}. The last part of the paper is devoted to highlight the impact of our findings in this research area. 

\vspace{0.3cm}

\end{abstract}

\section{Introduction}\label{Sect:Intro}


For fixed $n=2,3,\dots$ and $r=1,2,\dots$, we consider $r$ particles that are randomly distributed among $n$ cells  and let 
\[
X_j = \mbox{ number of particles fallen in cell $j$} \quad \mbox{ for } \quad j=1, \dots, n \,.
\]
The random variables $X_1, \dots, X_n$ are called \emph{occupancy numbers} and their joint distribution is said an \emph{occupancy model}. Note that the vector $\mathbf{X} \equiv \left( X_{1}, X_{2}, \dots, X_{n} \right)$ takes values in the set $A_{n,r}$, defined by
\begin{equation}\label{support:set}
A_{n,r} := \left\{ \mathbf{x} \equiv \left(  x_{1},\dots,x_{n} \right) : x_{j} \in \left\{ 0, 1, \dots, r \right\}  \mbox{ and } \sum_{j=1}^{n} x_{j}=r \right\} \,.
\end{equation}
An occupancy model is then a probability distribution on $A_{n,r}$ and we denote by $\mathcal{P} \left( A_{n,r} \right)$ the family of such probability distributions. Furthermore, the occupancy model is \emph{exchangeable} if $X_1, \dots, X_n$ are exchangeable; i.e., $X_1, \dots, X_n$ has the same distribution of $X_{\sigma(1)}, \dots, X_{\sigma(n)}$, whenever $\sigma$ is a permutation of the indices $\{1,2,\dots,n\}$. \\
%
%
The study of problems of placing objects in cells naturally arises in Statistical Physics in describing the microscopic arrangement of $r$ particles (which might be protons, electrons, \dots) among $n$ states (which might be energy levels). In this context there are three well-known examples of exchangeable occupancy distributions: Maxwell-Boltzmann (MB), Bose-Einstein (BE) and Fermi-Dirac (FD). 
One says that distinguishable (indistinguishable, respectively) particles that are not subject to the Pauli exclusion principle\footnote{At most one particle in each state.} obey MB statistics (BE, respectively). If, however, the particles are indistinguishable and subject to exclusion principle, they obey FD statistics. For example, electrons, protons and neutrons obey FD statistics; whereas photons and pions obey BE statistics.\\
In the sequel we will refer to MB, BE and FD as the {\em classical} occupancy models. {%
In Table~\ref{tab:classical_EOMs} of Section~\ref{Sect:Ma} 
}%
a more mathematical description of their main features is given. See \citet{Fel68} for further details.\\
It is interesting to understand why such models have a so fundamental role within the family of EOMs. In this respect, as an instrument to investigate and characterize what are the features that may justify the privileged position of MB, BE and FD statistics, in \cite{CoLeSpSu13} the authors introduced a remarkable subclass of EOMs: the $\mathcal{M}^{(a)}_{n,r}$-models.\\
{%
The set of $\mathcal{M}^{(a)}_{n,r}$-models is a wide family of EOMs, where every element of the class  can be obtained through conditionally i.i.d. random variables with laws in the exponential family. 
Besides, in particular, it contains the classical models. \\
We point out that this family is of interest in its own. 
\cite{CoLeSpSu13} focus on EOMs and their relations with the Uniform Order Statistics Property (UOSP) for discrete time counting processes. 
They show how definitions and results presented  in \citet{HuSh94}  and \citet{ShSpSu04} can be unified and generalized in the frame of $\mathcal{M}^{(a)}_{n,r}$-models. 
In particular, they allow to define a \emph{generalized UOSP} for the conditional distribution of the jump amounts given a certain number of arrivals. 
For processes with this property, they prove several characterizations in terms of $\mathcal{M}^{(a)}_{n,r}$-type distributions.\\
Moreover, the authors show that the $\mathcal{M}^{(a)}$-class is closed under some natural transformations of EOMs as the dropping of particles and the conditioning on a given occupancy sub-model. 
See Section \ref{Sect:Merging} for precise definitions and related results. \\
In this paper we still concentrate our attention on $\mathcal{M}^{(a)}_{n,r}$-models. 
More precisely, we introduce and study a new transformation of occupancy models, called \emph{merging}, that turns out to play a major role within the family of $\mathcal{M}^{(a)}_{n,r}$-models. 
Our aim is twofold. 
On the one hand, we investigate the algebraic properties of such application and how they reflect on the $\mathcal{M}^{(a)}$-class. 
On the other, we exploit these results to explain the relation between occupancy distributions of the type $\mathcal{M}^{(a)}$ and scale-consistent maximum entropy problems. \\
From the probabilistic point of view, the analysis of the merging transformation in the context of $\mathcal{M}^{(a)}_{n,r}$-models leads to a better understanding of the structure of this class. 
In fact, it suggests a decomposition in two subclasses: the first comprises {\it indecomposable} models, i.e. models that cannot be obtained through merging; the second consists of  the remaining ones, that is \emph{decomposable} models, derivable by merging. 
For instance, BE and FD distributions give rise to indecomposable occupancy models; whereas, a well known example of decomposable distribution is the {\it pseudo-contagious} model introduced by \citet{Cha05} and reconsidered by \cite{CoLeSpSu13}. 
Other examples are given in Section \ref{Sect:Merging}. \\
The decomposition just described above relies on the closure of the class $\mathcal{M}^{(a)}$ under merging and moreover, on the complete characterization we can provide of the occupancy model obtained after applying this transformation.\\
As previously mentioned, these are relevant facts for our purposes, since they are of potential interest in applications. Indeed, we conclude the manuscript  by showing a connection between $\mathcal{M}^{(a)}_{n,r}$-models and {\it entropy maximization (MaxEnt) inference}. \\
\citet{HaEt10} observe that the solution of MaxEnt depends on the scale at which the method is applied and, as a consequence, highlight the relevance of the choice of the observation scale when estimating an occupancy distribution. 
For this reason, they introduce the concept of \emph{scale-consistency}.
Roughly speaking, a MaxEnt distribution is self-consistent in the case when, if a different coarse-graining on cell size is considered, the type of probability distribution that solves the maximization problem does not change.\\
\citet{HaEt10} obtain MB and BE as solutions of some MaxEnt problems which are relevant in ecology for estimating the abundance of species. Moreover, as far as scale-consistency is concerned, they show that MB satisfies this property and BE does not.\\
A crucial remark is that our merging transformation coincides with the change of coarse-graining on cell size (\emph{scale transformation}) suggested by \citet{HaEt10}.
Therefore, in view of this result, we make clear (and formally prove) the reason why MB is scale-consistent, while BE is not. Furthermore, we are able to exhibit a whole family of self-consistent distributions. See Section~\ref{Sect:Merging} and Section~\ref{Sect:MaxEnt}, respectively.\\

}%

The outline of the paper is as follows. 
In Section~\ref{Sect:Ma} we recall the class of $\mathcal{M}^{(a)}$-models along with some of its basic features {and some related results in the framework} of discrete time counting processes. Section~\ref{Sect:Merging} is devoted to defining the merging transformation and to presenting 
our results concerned with its characterization. 
In addition, we also analyze the interplay between the new proposal and transformations previously introduced in \cite{CoLeSpSu13}. 
In Section~\ref{Sect:MaxEnt} we discuss the role of $\mathcal{M}^{(a)}_{n,r}$-models in the  { context} of entropy maximization inference. 
All proofs, if not immediate, are collected in Appendix~\ref{Sect:Proofs}. 

\section{The family $\mathcal{M}^{(a)}$}\label{Sect:Ma}

This family is parametrized by a function \mbox{$a: \mathbb{N} \longrightarrow ]0,+\infty[$} and $\left( X_1, \dots, X_n \right)$ is distributed according to the exchangeable occupancy model $\mathcal{M}_{n,r}^{(a)}$ if its distribution $P \in \mathcal{P} \left( A_{n,r}\right)$ has the form 
\begin{equation}\label{def:M(a)-model}
P \{ X_1=x_1, \dots, X_n=x_n \} = \frac{\prod_{j=1}^n a(x_j)}{C_{n,r}^{(a)}} \,,
\end{equation}
where 
\begin{equation}\label{norm:const}
C_{n,r}^{(a)} = \sum_{\boldsymbol{\eta} \in A_{n,r}} \prod_{j=1}^n a \left( \eta_j \right)
\end{equation} 
is the normalizing constant.\\
We recall that the classical models can be recovered with specific choices of the function $a$. Notice that they are respectively obtained by letting
\begin{align}
&\mbox{MB:} & &a(x) = \dfrac{1}{x!}, & &\mbox{for } x=0, 1, 2, \dots \label{a:MB}\\
&\mbox{BE:}  &  &a(x) = 1, & &\mbox{for }  x = 0, 1, 2, \dots \label{a:BE}\\
&\mbox{FD:}  &  &a(0)=a(1)=1, \quad a(x)=0, & &\mbox{for } x = 2, 3, \dots \label{a:FD}
\end{align}

{ %
We summarize the main characteristics of such models in Table~\ref{tab:classical_EOMs} below. 
}%

\begin{table}[h!]
\centering
\begin{tabular}{|c|c|l|c|}
\hline
{\small \bf Model} & {\small \bf Allocation Procedure} & {\small \bf Occupancy Distribution} & {\small \bf Support} \\
\hline
MB & \parbox{6cm}{\small \begin{itemize} \item $n$ distinguishable cells with \emph{unlimited} capacity \\ \item $r$ \emph{distinguishable} particles \\ \item the particles are distributed uniformly at random among cells \end{itemize}}  & $\displaystyle{P \{ \mathbf{X}=\mathbf{x} \} = \binom{r}{x_1 \cdots x_n} \, \frac{1}{n^{r}}}$ & $\mathbf{x} \in A_{n,r}$\\
\hline
BE & \parbox{6cm}{\small \begin{itemize} \item $n$ distinguishable cells with \emph{unlimited} capacity \\ \item $r$ \emph{indistinguishable} particles \\ \item the particles are distributed uniformly at random among cells \end{itemize}} & $\displaystyle{P \{ \mathbf{X}=\mathbf{x} \} = \binom{n+r-1}{n-1}^{-1}}$ & $\mathbf{x} \in A_{n,r}$\\
\hline
FD & \parbox{6cm}{\small \begin{itemize} \item $n$ distinguishable cells with capacity of \emph{at most one} particle\\ \item $r$ \emph{indistinguishable} particles \\ \item the particles are distributed uniformly at random among cells \end{itemize}} & $\displaystyle{P \{ \mathbf{X}=\mathbf{x} \} = \binom{n}{r}^{-1}}$ & $\mathbf{x} \in \widehat{A}_{n,r}$\\
\hline
\end{tabular}
\caption{Main features of the classical occupancy models. The set $A_{n,r}$ is as defined in \eqref{support:set}, while  $\widehat{A}_{n,r}$ is the subset of $A_{n,r}$ such that $x_{j} \in \{ 0,1 \}$ for all $j=1,\dots,n$.}
\label{tab:classical_EOMs}
\end{table}


It is interesting to highlight the connection between the class of the $\mathcal{M}_{n,r}^{(a)}$-models and the exponential family of discrete probability distributions. We believe this relationship also makes clear the role of the function $a$.\\
Let $a: \mathbb{N} \longrightarrow ]0,+\infty[$ and $b: \mathbb{R}_{+} \longrightarrow ]0,+\infty[$ be two fixed functions. We consider a sequence $V_1, V_2, V_3, \dots$ of $\mathbb{N}$-valued random variables and we suppose they are i.i.d. conditionally on a positive real parameter $\theta$, with univariate conditional marginal of the form
\begin{equation}\label{exp_fam:uni_marg}
f (v \vert \theta) = b (\theta) a (v) e^{- \theta v} .
\end{equation}
Thus, we get the $n$-dimensional joint distribution
\begin{equation}\label{exp_fam:joint_distr}
f^{(n)} \left( v_{1}, \dots, v_{n} \right) = \int_0^{+\infty} \left[ b(\theta) \right]^n \prod_{j=1}^{n} a \left( v_{j} \right) e^{- \theta \sum_{j=1}^n v_{j}} \, \Lambda (d\theta) \,,
\end{equation}
with $\Lambda$ a probability density on the positive real half-line. \\
By means of the variables $V$'s we may obtain a $\mathcal{M}_{n,r}^{(a)}$-model \eqref{def:M(a)-model} as follows. Set $S_n = \sum_{j=1}^n V_{j}$ and consider the random vector $(X_1, \dots, X_n)$ with values in $A_{n,r}$, whose distribution is given by
\[
P \left\{ X_1 = x_1, \dots, X_n = x_n \right\} := P \left\{ \left. V_{1}=x_1, \dots, V_{n}=x_n \right\vert S_n= r\right\} \quad \mbox{ for } \mathbf{x} \in A_{n,r} \,.
\]
Thanks to \eqref{exp_fam:joint_distr} it is readily seen that
\[
P \left\{ X_1 = x_1, \dots, X_n = x_n \right\} = \frac{\prod_{j=1}^n a (x_j)}{C_{n,r}^{(a)}}  \,,
\]
for any distribution $\Lambda$. This example suggests a few observations. 

\begin{remark}
If the distribution $\Lambda$ is degenerate, the construction of $\mathcal{M}_{n,r}^{(a)}$-models through i.i.d. random variables in \citet{Cha05} is recovered. 
\end{remark}

\begin{remark}
The three classical models can be obtained by choosing suitably the distribution $f(v \vert \theta)$. In fact, MB, BE and FD are generated by a Poisson, Geometric and Bernoulli distribution, respectively.
\end{remark}

\begin{remark}
The fact that $P \{ V_{1} = x_{1}, \dots, V_{n} = x_{n} | S_{n} = r \}$ does not depend on the distribution $\Lambda$, has an immediate interpretation in statistical terms: for the \emph{exponential model}, $S_{n}$ is a sufficient statistic with respect to the parameter $\theta$. 
\end{remark}

\begin{remark}
Up to a multiplicative factor, the normalization constant $C_{n,r}^{(a)}$ can be interpreted as the probability that a sum  of certain i.i.d. random variables is equal to $r$, i.e. $C_{n,r}^{(a)} = \kappa P\{S_n=r\}$ with $\kappa$ suitable constant. 
\end{remark}


\bigskip

$\mathcal{M}^{(a)}_{n,r}$-models play a major role in the context of discrete time counting processes. They provide a framework where defining a natural and unifying extension of the Uniform Order Statistics Property (UOSP) given in \citet{HuSh94} and the UOSP($\leq$) introduced by \citet{ShSpSu04}. Precisely, the generalization proposed by \citet{CoLeSpSu13} and called $\mathcal{M}^{(a)}$-UOSP is as follows.

\begin{definition}\label{Def:Ma-UOSP}
Let $a:\{0,1,\dots\} \longrightarrow ]0,+\infty[$ be a given function and $\{ N_{t} \}_{t = 0, 1, \dots}$ a discrete-time counting process with jump amounts $J_0, \dots, J_t$. We say that it satisfies the $\mathcal{M}^{(a)}$-UOSP if, for any $t, k \in \mathbb{N}$ and any $(j_{0}, \dots, j_{t}) \in A_{t+1,k}$, we have
\begin{equation}\label{Def:M-UOSP}
P \{ J_{0} = j_{0}, \dots, J_{t} = j_{t} | N_{t} = k \} = \frac{\prod_{h=0}^t a(j_{h})}{C_{t+1,k}^{(a)}}.
\end{equation}
\end{definition}  

It is readily seen that, by choosing function $a$ in Definition~\ref{Def:Ma-UOSP} as \eqref{a:FD} and \eqref{a:MB}, we can respectively recover the UOSP and UOSP($\leq$).\\
In a similar fashion of \citet{HuSh94} and \citet{ShSpSu04}, for processes satisfying the $\mathcal{M}^{(a)}$-UOSP, several characterizations may be proven. To state exhaustively these results, we need to recall first the notion of $a$-mixed geometric process, which, in turn, can be seen as an appropriate generalization of the definition of mixed geometric process in \citet{HuSh94}.

\begin{definition}
Let $a: \{ 0, 1, \dots \} \longrightarrow ]0,+\infty[$ be a given function such that $a(0)=1$. The process $\{ N_{t} \}_{t = 0, 1, \dots, M}$ is an $a$\emph{-mixed geometric process} if the discrete joint density of $(J_{0}, J_{1}, \dots, J_{t} )$ has the form
\begin{equation}\label{CondDefMixedaProc}
p_{t} (j_{0}, j_{1}, \dots, j_{t}) = R_{t} \left( \sum_{h=0}^{t} j_{h} \right) \cdot \prod_{h=0}^{t} a(j_{h}) \quad \mbox{for } t = 0, 1, \dots, M
\end{equation}
for a suitable sequence of functions $R_{t}: \{0, 1, \dots \} \longrightarrow \mathbb{R}_{+}$, $t=0,1,\dots,M$.
\end{definition} 

On the one hand, \citet{CoLeSpSu13} prove  the equivalence for a multiple jumps counting process of being an $a$-mixed geometric process and satisfying the $\mathcal{M}^{(a)}$-UOSP; on the other hand, for processes fulfilling \eqref{Def:M-UOSP} (or, equivalently, \eqref{CondDefMixedaProc}) they are able to give a characterization of the joint distribution of both arrival and inter-arrival times.\\

Anyway, this is not the only application of $\mathcal{M}^{(a)}_{n,r}$-models. For instance, in Section \ref{Sect:MaxEnt} we will highlight the importance of occupancy models in the so-called entropy maximization inference. The results provided in the next Section~\ref{Sect:Merging} are crucial to explain some properties in MaxEnt framework from a probabilistic point of view. 

\section{The merging transformation}\label{Sect:Merging}

Fix $N = 2, 3, \dots$ and $r = 1, 2, \dots$ and let $s$ be a proper divisor of $N$, i.e. $N = ns$ with $n, s \in \mathbb{N}$, $s \neq N$. Consider an occupancy model over $A_{N,r}=A_{ns,r}$ and denote by 
\[
\mathbf{Z}_1^{(s)} := ( Z_{1,1}, \dots, Z_{1,s}) \quad \mathbf{Z}_2^{(s)}  := (Z_{2,1}, \dots, Z_{2,s}) \quad \dots \quad \mathbf{Z}_n^{(s)} := (Z_{n,1}, \dots, Z_{n,s})
\]
the corresponding occupancy numbers.
 We merge groups of $s$ cells to create $n$ macrocells. We obtain a new  occupancy model on $A_{n,r}$, whose occupancy numbers are $X_1, \dots, X_n$ with  
\[
X_j = \sum_{h=1}^s Z_{j,h} 
\]
and joint distribution
\[
P \{ X_1=x_1, \dots, X_n=x_n \} = \sum_{\mathbf{z} \in H_{n,r, \mathbf{x}}^{(s)}} P\left\{ Z_{1,1}=z_{1,1}, \dots, Z_{n,s}=z_{n,s} \right\} \,,
\]
where
\begin{equation}\label{Constraints:set}
H_{n,r, \mathbf{x}}^{(s)} := \left\{ \mathbf{z} \equiv \left( \mathbf{z}_1^{(s)}, \dots, \mathbf{z}_n^{(s)} \right) \in A_{ns,r} : \mathbf{z}_j^{(s)} \in A_{s,x_j} \mbox{ for } j=1, \dots, n  \mbox{ and given } \mathbf{x} \in A_{n,r}  \right\} \,.
\end{equation}
We denote the merging operation described above with the symbol $I_s$. Let us now turn to considering the case of exchangeability. Notice, first of all, that such a condition is actually preserved under the operation $I_{s}$. It is remarkable, furthermore, that the same EOM is obtained even if we form the n macrocells by   choosing in a whatever different way the n groups of cells to be merged, provided any group still contains exactly s cells. A further invariance property is described by the following result.

\begin{proposition}\label{Prop:composition:I}
Let $P \in \mathcal{P} \left( A_{N,r} \right)$ be an EOM. Fix $s_1, s_2$ proper divisors of $N$, such that $s_1$ and $s_2$ are coprime ($N = n s_1 s_2$); then, 
\begin{equation}\label{composition:I}
\left( I_{s_1} \circ I_{s_2} \right) (P) = I_{s_1s_2} (P) \,.
\end{equation}
\end{proposition}

\begin{remark}
The divisors $s_1$ and $s_2$ are taken coprime to ensure that $N$ can be sequentially divided by them.
\end{remark}

By induction it is easy to see that the result in Proposition \ref{Prop:composition:I} can be extended up to a sequence $s_1, \dots s_k$ of mutually coprime proper divisors of $N$. Therefore, the representative elements in the class of exchangeable occupancy models are those with a prime number of cells, since all other models within this subset may be iteratively merged up to reach a model on $A_{n,r}$ with $n$ prime. \\

Now we focus on the class $\mathcal{M}^{(a)}$. It is natural to wonder what occurs whenever we apply transformation $I_s$ to an occupancy model of the form \eqref{def:M(a)-model}. The following proposition answers this question. On the one hand, it proves that the family $\mathcal{M}^{(a)}$ is closed under $I_s$. On the other, it shows that it is possible to exploit the special structure of the probability distribution \eqref{def:M(a)-model} to completely characterize the resulting merged occupancy model.

\begin{proposition}\label{Prop:closure:I}
Let $P \in \mathcal{P} \left( A_{N,r} \right)$ be a $\mathcal{M}^{(a)}_{N,r}$-model. Fix $s$ a proper divisor of $N$ ($N=ns$); then, $P' = I_s (P) \in \mathcal{P} \left( A_{n,r} \right)$, obtained by applying the transformation $I_s$, is distributed according to a model $\mathcal{M}^{(a')}_{n,r}$, with $a' (x) = C^{(a)}_{s,x}$.
\end{proposition}

To better explain the consequences of previous propositions, let us fix some notation. We denote by $\mathcal{I}_{\mathcal{M}^{(a)}}$ the set of occupancy models obtained as follows: 
\begin{equation}\label{def:class:IMa}
\framebox{\parbox{0.6\textwidth}{
The elements of $\mathcal{I}_{\mathcal{M}^{(a)}}$ are those $\mathcal{M}^{(a)}$-models of the form \\[0.1cm] $P'=I_s (P)$ for some $s$ and some $P$, with $P$ of the type $\mathcal{M}^{(a)}$.
}}
\end{equation}

From Proposition~\ref{Prop:composition:I} and Proposition~\ref{Prop:closure:I} we can deduce that the class $\mathcal{I}_{\mathcal{M}^{(a)}}$ is closed with respect to the merging transformation. Let $s_1, s_2$ be proper divisors of $N$, such that $s_1$ and $s_2$ are coprime ($N = n s_1 s_2$), and consider $P$ an occupancy model in the class $\mathcal{M}^{(a)}_{N,r}$. We have the following mapping
\[
\xymatrix@R=3pt@C=40pt{
\text{\parbox{3cm}{ \centering $P \in \mathcal{M}^{(a)}_{ns_1s_2,r}$-class \\ with \\ $a(x)$}} \ar@{->}@//[r]^-{I_{s_1}} & \text{\parbox{3.5cm}{ \centering $I_{s_1} (P) \in \mathcal{M}^{(a')}_{ns_2,r}$-class \\ with \\ $a'(x)=C^{(a)}_{s_1,x}$}} \ar@{->}@//[r]^-{I_{s_2}} & 
\text{\parbox{4.5cm}{ \centering $(I_{s_2} \circ I_{s_1}) (P) \in \mathcal{M}^{(a'')}_{n,r}$-class \\ with \\ $a''(x)=C^{(a')}_{s_2,x}=C^{(a)}_{s_1s_2,x}$}}
}
\]
The closure property arises from the fact that $a'$ and $a''$ have the same form as functions of $x$; however they depend on different scale parameters $s_1$ and $s_1s_2$, respectively. Therefore, if we iteratively apply the $I_s$ transformation, at the first step we \emph{possibly} modify the function $a$ characterizing the model, but then we will affect only the scale parameter. Roughly speaking, starting from a given occupancy $\mathcal{M}^{(a)}_{n,r}$-model, we are constructing a sort of cone within the class $\mathcal{M}^{(a)}$ through merging.  \\

We believe it is worth to illustrate explicitly the three classes $\mathcal{I}_{\mathrm{MB}}$, $\mathcal{I}_{\mathrm{BE}}$ and $\mathcal{I}_{\mathrm{FD}}$, obtained by merging the classical models. The formal definitions of these families are analogous to \eqref{def:class:IMa}.

\begin{description}
\item[Class $\mathcal{I}_{\mathrm{MB}}$] --- Consider a MB model on the set $A_{ns,r}$, that is
\[
P\left\{ Z_{1,1}=z_{1,1}, \dots, Z_{1,s}=z_{1,s}, \dots, Z_{n,1}=z_{n,1}, \dots, Z_{n,s}=z_{n,s} \right\} = \frac{1}{(ns)^r} \frac{r!}{z_{1,1}! \cdots z_{n,s}!} \,,
\]
and then apply transformation $I_s$. It readily yields
\[
P \{ X_1=x_1, \dots, X_n=x_n \} 
= \frac{1}{n^r} \frac{r!}{x_1! \cdots x_n!} \,.
\]
Notice that the merged configuration we obtained is still distributed as a MB occupancy model, but over the set $A_{n,r}$ now. It means the class of MB distributions is itself closed under transformation $I_s$ and $\mathcal{I}_{\mathrm{MB}}$ just coincides with the class of all MB models. \\

\item[Class $\mathcal{I}_{\mathrm{BE}}$] --- Consider a BE model on the set $A_{ns,r}$, that is
\[
P \left\{ Z_{1,1}=z_{1,1}, \dots, Z_{1,s}=z_{1,s}, \dots, Z_{n,1}=z_{n,1}, \dots, Z_{n,s}=z_{n,s} \right\} = \frac{1}{{ns+r-1 \choose ns-1}}\,,
\]
and apply transformation $I_s$. It readily yields
\begin{equation}\label{pseudo-contagious}
P \{ X_1=x_1, \dots, X_n=x_n \} 
= \frac{{s+x_1-1 \choose x_1}{s+x_2-1 \choose x_2} \cdots {s+x_n-1 \choose x_n}}{{ns+r-1 \choose r}} \,.
\end{equation}
In view of \eqref{pseudo-contagious}, we point out that the class $\mathcal{I}_{\mathrm{BE}}$ coincides with the class of \emph{pseudo-contagious} occupancy models presented in \citet{Cha05}. \\

\item[Class $\mathcal{I}_{\mathrm{FD}}$] --- Consider a FD model over the set $\widehat{A}_{ns,r} := \left\{ \mathbf{z} \in A_{ns,r} : z_j \in \{0,1\} \right\}$, that is
\[
P\left\{ Z_{1,1}=z_{1,1}, \dots, Z_{1,s}=z_{1,s}, \dots, Z_{n,1}=z_{n,1}, \dots, Z_{n,s}=z_{n,s} \right\} = \frac{1}{{ns \choose r}} \,,
\]
and apply transformation $I_s$. It readily yields
\[
P \{ X_1=x_1, \dots, X_n=x_n \} 
=\frac{{s \choose x_1}{s \choose x_2} \cdots {s \choose x_n}}{{ns \choose r}} \,, 
\]
with $ \left( x_1, \dots, x_n \right) \in \widetilde{A}_{n,r} := \left\{ \mathbf{x} \in A_{n,r} : x_j \in \{0, \dots, s \} \right\}$. The class $\mathcal{I}_{\mathrm{FD}}$ contains all \emph{multi-hypergeometric} occupancy models.
\end{description}


\begin{remark}\label{Rmk:closure:I:fundamental:classes}
From Proposition \ref{Prop:composition:I} and Proposition \ref{Prop:closure:I} follows that the classes $\mathcal{I}_{\mathrm{BE}}$, $\mathcal{I}_{\mathrm{FD}}$ are themselves closed under transformation $I_s$. 
\end{remark}

\begin{remark}
Observe that applying $I_s$ within the classes $\mathcal{I}_{\mathrm{BE}}$, $\mathcal{I}_{\mathrm{FD}}$ and $\mathcal{I}_{\mathcal{M}^{(a)}}$ means observing a given occupancy model on a different scale. 
\end{remark}


From a heuristic point of view, the class $\mathcal{M}^{(a)}$ is comprised of models which are a sort of generators, in the sense that they cannot be constructed as merging of other elements in the class --- we naively call them \emph{germs} ---, and the models obtained from those germs through merging (subset $\mathcal{I}_{\mathcal{M}^{(a)}}$).  We may look at the family $\mathcal{M}^{(a)}$ as follows 
\begin{equation}\label{Ma:decomposition}
\mathcal{M}^{(a)} = \mathcal{G}_{\mathcal{M}^{(a)}} \cup \mathcal{I}_{\mathcal{M}^{(a)}} \,,
\end{equation}
where $\mathcal{G}_{\mathcal{M}^{(a)}}$ is the set of all germs. From Proposition~\ref{Prop:closure:I} the following circumstance emerges:  in order to understand if a model in $\mathcal{M}^{(a)}$, with a given $a$, is a \emph{germ} or a \emph{merged model} it suffices to analyze the function $a$. Indeed, a considered $\mathcal{M}^{(a)}$-model is a germ if and only if the corresponding $a$ is not a power of convolution; in other words, if and only if $a$ gives rise to an \emph{indecomposable distribution} (see \cite{LiOs77}). 

\begin{remark}
Note that $\mathcal{I}_{\mathcal{M}^{(a)}} \subsetneq \mathcal{M}^{(a)}$.  For instance, BE and FD models cannot be obtained from other elements in the $\mathcal{M}^{(a)}$ class through merging. 
\end{remark}

\begin{figure}[h]
\centering%
\includegraphics[width=.6\textwidth]{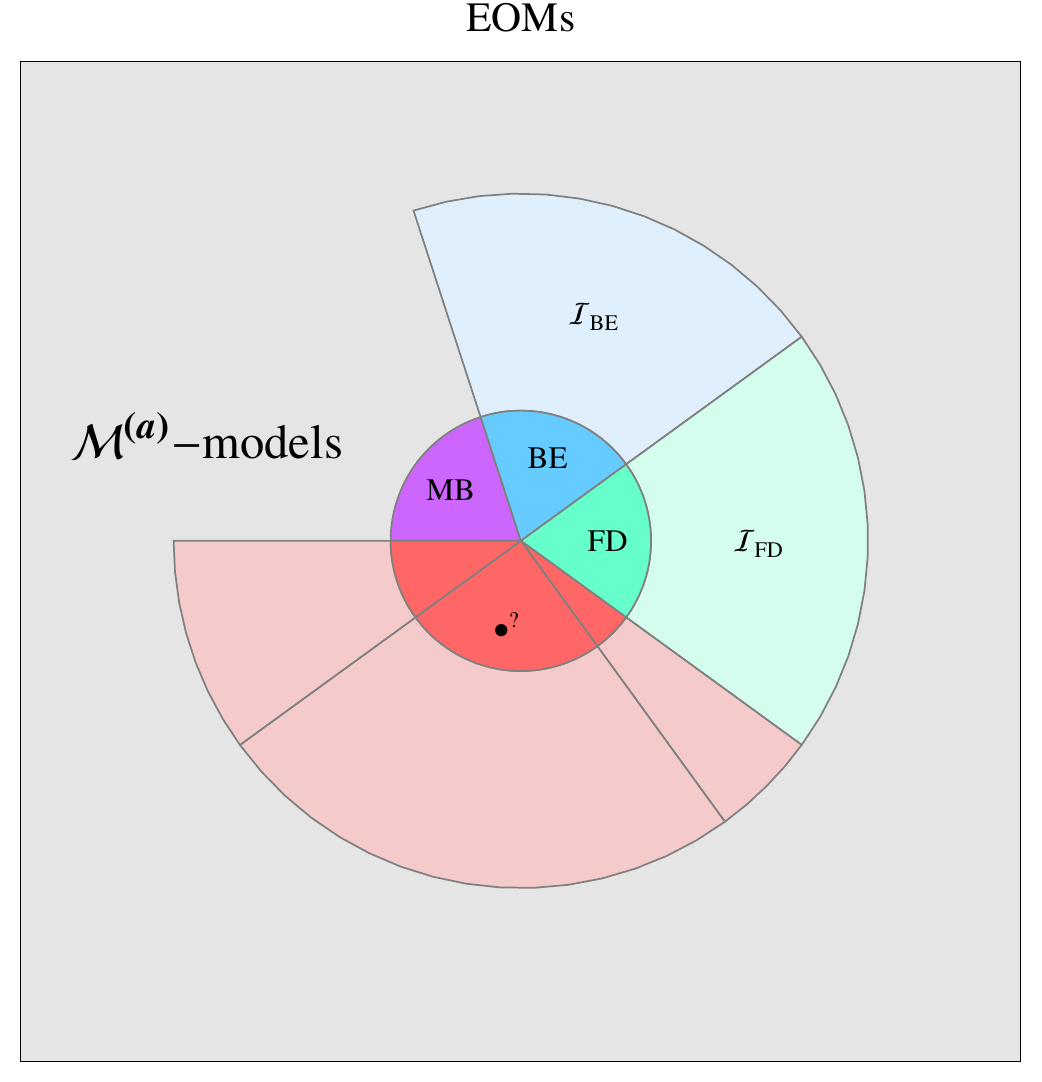}
\caption{Schematic view of the $\mathcal{M}^{(a)}$-models within the EOMs. The set of regions in solid color are the germs, whereas the dissolved ones altogether form the class $\mathcal{I}_{\mathcal{M}^{(a)}}$. Each sector represents the cone generated by a germ through the application of transformation $I_s$.}
\label{fig:1}
\end{figure}

In Figure~\ref{fig:1} an illustration of the structure of the $\mathcal{M}^{(a)}$-class is given. 

It is then natural to wonder if there exist further germs, different from MB, BE and FD, that generate other cones within the $\mathcal{M}^{(a)}$ class through transformation $I_s$. Due to the arbitrariness of the function $a$ however, we immediately see there are infinitely many germs. \\
In general, for a particular given $a$, it is not trivial to characterize the elements of the corresponding set $\mathcal{I}_{\mathcal{M}^{(a)}}$. To be concrete let us mention the following example.

\begin{example}
Take the model corresponding to the function $a(x) = \widehat{a}(x)=\frac{1}{x+1}$ on $A_{ns,r}$. Observe that it may be constructed starting from conditionally i.i.d. random variables $V_{1,1}, \dots, V_{n,s}$ by choosing as \eqref{exp_fam:uni_marg} the logarithmic series distribution (see \citet{JoKo69} for more details)
\[
f(v \vert \theta) = - \frac{\theta^{v+1}}{(v+1)\log(1-\theta)} \quad \mbox{ with } \quad 0 < \theta < 1 \,.
\] 
Now apply transformation $I_s$. We know (see Proposition~\ref{Prop:closure:I}) that the resulting occupancy model belongs to the class $\mathcal{M}_{n,r}^{(\widehat{a}')}$, with $\widehat{a}'(x) = C_{s,x}^{(\widehat{a})}$. In particular, for a suitable positive constant $\kappa$, one has
\[
C_{s,x}^{(\widehat{a})} = \kappa \, P \left\{ \sum_{h=1}^{s} V_{j,h} = x\right\} \,.
\]
Notice however that we are not able to give an explicit form to the probability in the right-hand side.\\
\end{example}

The above example enlightens the following fact: in order to describe the function $a'$ characterizing a merged model, it is crucial being able to determine the distribution of a sum of certain (conditionally) i.i.d. random variables.\\ 
The interest in the classes MB, BE, FD, $\mathcal{I}_{\mathrm{BE}}$ and $\mathcal{I}_{\mathrm{FD}}$ is in the fact that the elements in there are those generated by distributions for which it is possible to determine precisely such quantity. Table~\ref{tab:convolution} below makes evident the reason why the MB, $\mathcal{I}_{\mathrm{BE}}$ and $\mathcal{I}_{\mathrm{FD}}$ classes are closed under merging. 

\begin{table}[h!]
\centering
\begin{tabular}{|c|c|c|c|c|}
\hline
{\small \bf Class} & {\small \bf Underlying Distribution} &  {\small \bf $n$-convolution} \\
\hline
&&\\
MB & $f(v \vert \lambda) \sim \mathrm{Po} (\lambda)$ & $f^{*n}(v \vert \lambda) \sim \mathrm{Po} (n\lambda)$ \\
&&\\
\hline
&&\\
BE & $f(v \vert p) \sim \mathrm{Ge} (p)$ & $f^{*n}(v \vert p) \sim \mathrm{NB} (n,p)$ \\
&&\\
\hline
&&\\
FD & $f(v \vert p) \sim \mathrm{Be} (p)$ & $f^{*n}(v \vert p) \sim \mathrm{B} (n,p)$ \\
&&\\
\hline
&&\\
$\mathcal{I}_{\mathrm{BE}}$ & $f(v \vert p) \sim \mathrm{NB} (m,p)$ & $f^{*n}(v \vert p) \sim \mathrm{NB} (nm,p)$ \\
&&\\
\hline
&&\\
$\mathcal{I}_{\mathrm{FD}}$ & $f(v \vert p) \sim \mathrm{B} (m,p)$ & $f^{*n}(v \vert p) \sim \mathrm{B} (nm,p)$\\
&&\\
\hline
\end{tabular}
\label{tab:convolution}
\caption{Convolutions of the underlying distribution for the classical models and the corresponding merged classes.}
\end{table}

\subsection{Relations between $I_s$ and other classes of transformations}

To understand the structure of classes of EOMs, an approach is to investigate closure properties under some natural transformations. We briefly recall transformations and related results presented in \citet{CoLeSpSu13}.

\begin{description}
\item[Transformation $K_{1}: \mathcal{P} \left( A_{n,r} \right) \longrightarrow \mathcal{P} \left( A_{n,r-1} \right)$.] Consider $r$ particles distributed among $n$ cells according to a given occupancy model. We drop one of the particles from this population uniformly at random and we obtain a new occupancy model with $n$ cells and $r-1$ particles. \\

\item[Transformation $K_{n,r}^{N,R}: \mathcal{P} \left( A_{N,R} \right) \longrightarrow \mathcal{P} \left( A_{n,r} \right)$.] Consider $R$ particles distributed among $N$ cells according to a given occupancy model and let $X_1, \dots, X_N$ be the related occupancy numbers. Here and in the rest we will use the notation $S_N = \sum_{j=1}^N X_j$. Fix positive integers $n$ and $r$ such that $1 \leq n \leq N-1$, $1 \leq r \leq R$. We condition on the event $\{S_{n}=r\}$ and we obtain a new occupancy model with $n$ cells and $r$ particles, defined by
\[
P \{ X_{1} = x_{1}, \dots, X_{n} = x_{n} | S_{n} = r \} = \frac{P \{ X_{1} = x_{1}, \dots, X_{n} = x_{n} \}}{P \{ S_{n} = r \}} \quad \mbox{ for } \mathbf{x} \in A_{n,r}
\]
and where
\begin{multline*}
P \{ X_{1} = x_{1}, \dots, X_{n} = x_{n} \} \\
= \sum_{\boldsymbol{\eta} \in A_{N-n,R-r}} P \{ X_{1} = x_{1}, \dots, X_{n} = x_{n}, X_{n+1} = \eta_{1}, \dots, X_{N} = \eta_{N-n} \},
\end{multline*}
\[
P \{ S_{n} = r\} = \sum_{\mathbf{x} \in A_{n,r}} \sum_{\boldsymbol{\eta} \in A_{N-n,R-r}} P \{ X_{1} = x_{1}, \dots, X_{n} = x_{n}, X_{n+1} = \eta_{1}, \dots, X_{n} = \eta_{N-n} \}.
\]
\end{description}

Notice that all the considered transformations ($I_s$ included) are mappings from a set $\mathcal{P} \left( A_{n,r} \right)$ to a set $\mathcal{P} \left( A_{n',r'} \right)$, with $n' \leq n$ and $r' \leq r$ and preserve exchangeability. \\
It is relevant to look at what happens to the structure of $\mathcal{M}^{(a)}$-models whenever we apply one of these transformations. Indeed, the class of $\mathcal{M}^{(a)}$-models is closed under $K^{N,R}_{n,r}$; whereas, in the case we apply $K_1$ the structure of such class is preserved only if a technical condition is satisfied. 

\begin{proposition}[Proposition 6.5 in \citet{CoLeSpSu13}]\label{Prop:6.5}
Let $(X_1,\dots,X_N)$ be a $\mathcal{M}_{N,R}^{(a)}$-model and define $S_n=X_1+\cdots+X_n$ with $n\leq N$. Conditionally on the event $\lbrace S_n=r\rbrace$, the variables $X_1,\cdots, X_n$ are distributed as a $\mathcal{M}_{n,r}^{(a)}$-model.
\end{proposition}

\begin{proposition}[Proposition 6.6 in \citet{CoLeSpSu13}]\label{Prop:6.6}
Let $(X_{1}, \dots, X_{n})$ be a $\mathcal{M}_{n,r}^{(a)}$-model and let the function $a$ fulfill the condition
\begin{equation}\label{23}\tag{$\star$}
\frac{C_{n,r-1}^{(a)}}{C_{n,r}^{(a)}}\sum_{h=1}^{n} \frac{x'_{h}+1}{r}\frac{a(x'_{h}+1)}{a(x'_h)}=1 \,, \quad \mbox{ for any } {\bf x'} \in A_{n,r-1} \,.
\end{equation}
Then, the occupancy vector $(X_1',\dots,X'_n)$, obtained by applying the transformation $K_1$, is distributed according to the model $\mathcal{M}_{n,r-1}^{(a)}$.
\end{proposition}

Let us now turn to the class $\mathcal{I}_{\mathcal{M}^{(a)}}$, defined in \eqref{def:class:IMa}, and investigate possible closure properties under transformations $K_1$ and $K^{N,R}_{n,r}$. 

\begin{proposition}\label{Prop:dropping:I}
Consider a $\mathcal{M}^{(a')}_{n,r}$-model $P'$ such that $P' = I_s \left( P \right)$, for some $s$ and some $\mathcal{M}^{(a)}_{ns,r}$-model $P$. If the function $a'$ satisfies \eqref{23}, then $K_1 \left( P' \right) = I_s (\widetilde{P})$, for some $\mathcal{M}^{(a')}_{ns,r-1}$-model $\widetilde{P}$. \\
Besides, if in addition also function $a$ fulfills \eqref{23}, then 
\[
\left( K_1 \circ I_s \right) (P) = \left( I_s \circ K_1 \right) (P).
\]
\end{proposition}

The first statement in the previous proposition is concerned with the closure of the $\mathcal{I}_{\mathcal{M}^{(a)}}$ class with respect to transformation $K_1$. As corollaries we get some further results focusing on the special subclasses $\mathcal{I}_{\mathrm{BE}}$, $\mathcal{I}_{\mathrm{FD}}$ and MB. 

\begin{corollary}\label{Cor:dropping:I:subclasses}
Consider a $\mathcal{M}^{(a')}_{n,r}$-model $P'$ such that $P' \in \mathcal{I}_{\mathrm{BE}}$. Then $K_1 \left( P' \right)$ is a $\mathcal{M}^{(a')}_{n,r-1}$-model still belonging to the $\mathcal{I}_{BE}$ class.
\end{corollary}

In view of Proposition~\ref{Prop:dropping:I}, to prove Corollary~\ref{Cor:dropping:I:subclasses} it suffices to show that the function $a'(x) = {s+x-1 \choose x}$, corresponding to distribution \eqref{pseudo-contagious}, fulfills assumption \eqref{23} for every $s$. With notation introduced in Proposition~\ref{Prop:6.6}, we have
\[
\frac{C_{n,r-1}^{(a')}}{C_{n,r}^{(a')}}\sum_{h=1}^{n} \frac{x'_{h}+1}{r}\frac{a'(x'_{h}+1)}{a'(x'_h)} = \frac{1}{ns+r-1} \left( \sum_{h=1}^n s + x_h' \right) = 1 \quad \mbox{ for any } \mathbf{x'} \in A_{n,r-1},
\]
and this together with the assertion in Remark~\ref{Rmk:closure:I:fundamental:classes} leads to the conclusion. \\
In a similar fashion it is possible to check a result, analogous to Corollary~\ref{Cor:dropping:I:subclasses}, for models in the classes  $\mathcal{I}_{\mathrm{FD}}$ and MB.
%

Next result describes the interplay between transformations $I_s$ and $K^{N,R}_{n,r}$.

\begin{proposition}\label{Prop:conditioning:I}
Let $P$ a $\mathcal{M}^{(a)}_{Ns,R}$-model. Fix $n \leq N$, $r \leq R$ and $s$; then 
\begin{equation}\label{comm:cond:I}
\left( K_{n,r}^{N,R} \circ I_s \right) (P) = \left( I_s \circ K_{ns,r}^{Ns,R} \right) (P)
\end{equation}
and it is distributed as a $\mathcal{M}^{(a')}_{n,r}$-model.
\end{proposition}

\begin{remark}
Notice that all integer parameters entering in the transformations $K$'s in \eqref{comm:cond:I} are chosen so that both compositions are well-defined. 
\end{remark}

Previous proposition states in particular that the subset $\mathcal{I}_{\mathcal{M}^{(a)}}$ is closed under transformation $K^{N,R}_{n,r}$. Indeed, if $P'$ is a $\mathcal{M}^{(a')}_{N,R}$-model such that $P'=I_s(P)$, then identity \eqref{comm:cond:I} implies $K_{n,r}^{N,R} (P') = I_s (\widetilde{P})$, for some $\widetilde{P}$ in the set of $\mathcal{M}_{ns,r}^{(a)}$-models.
By taking into account this result and the arguments in Remark \ref{Rmk:closure:I:fundamental:classes},  as a corollary we obtain that the same closure property is true for the subclasses MB, $\mathcal{I}_{\mathrm{BE}}$ and $\mathcal{I}_{\mathrm{FD}}$. 

\section{Entropy maximization inference and $\mathcal{M}_{n,r}^{(a)}$-models}\label{Sect:MaxEnt}

{ %
The maximum entropy principle (MaxEnt) is an inference procedure, originated in statistical mechanics.  It allows one to determine the probability distribution of a microscopic unit configuration from macroscopic (observable and measurable) averaged quantities of a physical system. The core of the procedure relies on the circumstance that the probability distribution with highest entropy conveys the minimum information and best corresponds to the current state of knowledge about the system. Indeed, entropy provides a measure of uncertainty that has to be maximized when looking for the probability law that encodes only the available information. The latter is  represented by a set of constraints that the distribution itself must satisfy.  Typically such constraints are expected values of suitable functions of microscopic configurations. It is worth mentioning that, typically, they are exact averages with respect to the unknown probability distribution. \\
}%
See the seminal paper \cite{Jay57} for motivations and further details.\\
Let us be more precise and formalize MaxEnt in the discrete case. Let $S$ be a countable set of configurations and $p$ denote a probability distribution over $S$. Moreover, assume that the available information about the system can be cast in form of linear constraints with respect to probability (typically expected values). Solving MaxEnt amounts to solving the following optimization problem
\begin{align}\label{MaxEnt:pb}
\max_{p} & \quad - \sum_{s \in S} p(s) \log p(s) \nonumber\\
         & \mbox{subject to } \quad p(s) \geq 0 \quad \mbox{ for every } s \in S \nonumber\\
         & \phantom{subject to } \quad \sum_{s \in S} p(s) = 1 \\
         & \phantom{subject to } \quad \sum_{s \in S} f_k (s) p(s) = c_k \quad \mbox{ for } k=1, \dots, m \nonumber
\end{align}
where $f_k$'s and $c_k$'s are given real functions and constants, respectively. 

\begin{proposition}\label{Prop:MaxEnt:sol}
For fixed functions $f_k: S \longrightarrow \mathbb{R}$ and fixed values $c_k \in \mathbb{R}$, with $k=1, \dots, m$, the maximization entropy problem \eqref{MaxEnt:pb} admits the unique solution
\begin{equation}\label{Gibbs:distr:gen}
p^*(s) = \frac{e^{- \sum_{k=1}^m \lambda_k f_k(s)}}{\sum_{s \in S} e^{- \sum_{k=1}^m \lambda_k f_k(s)}} \,,
\end{equation}
where $\lambda_1, \dots, \lambda_m$ are suitable real constants that must be chosen so that $p^*$ fulfills the constraints. 
\end{proposition}

From expression \eqref{Gibbs:distr:gen} we immediately see that 
we can obtain different distributions 
by modifying constraints and support set. \\
Now let us set $S = A_{n,r}$, defined by \eqref{support:set}, and focus on occupancy distributions. 
We know that MB and BE (resp. FD) statistics are the ones maximizing the entropy in the cases when we are considering respectively probability distributions over $n$-tuples of distinguishable/indistinguishable (resp. indistinguishable, subject to exclusion principle) elements, see \cite{HaEt10}.\\  
Observe that also probability distributions of the form \eqref{def:M(a)-model} can be obtained as solutions of MaxEnt, provided appropriate functions $f_k$'s are selected. More precisely we can state

\begin{proposition}\label{Prop:MaxEnt:Ma}
The exchangeable occupancy model in \eqref{def:M(a)-model} is a solution of MaxEnt \eqref{MaxEnt:pb} over \mbox{$S = A_{n,r}$} under the constraints
\begin{align}\label{MaxEnt:Ma:constraints}
& p(\mathbf{x}) \geq 0 \quad \mbox{ for every } \mathbf{x} \in A_{n,r} \nonumber\\
& \sum_{\mathbf{x} \in A_{n,r}} p(\mathbf{x}) = 1 \\
& \sum_{\mathbf{x} \in A_{n,r}} \log \left( \prod_{j=1}^n a(x_j) \right) p(\mathbf{x}) = c \nonumber
\end{align}
where $a: \mathbb{N} \longrightarrow ]0,+\infty[$ is a fixed function and $c$ a suitable real constant. 
\end{proposition}

EOMs naturally arise in the framework of entropy maximization inference applied to ecological systems. For instance, imagine to investigate the spatial abundance distribution of a given species in a given region. We have at hand a huge amount of data consisting in positions of individuals in that region and we aim at unveiling the main features of the underlying structure of the data set. We want  to estimate the probability distribution from which the database is sampled.\\ 
We divide the considered geographical area in $n$ subregions (cells) and define a system configuration as the vector of occupancy numbers representing how many individuals (particles) of a given species are present in each subregion.  We then determine the  distribution of possible arrangements by maximizing the entropy under the natural constraint of having a population of specific size $r$.\\
Generally, when analyzing a dataset, a crucial issue is choosing the most informative \emph{scale} of a physical system. In particular, when dividing a geographical area into different cells, the appropriate choice of the number of such cells is essential for a correct inference about patterns in the data. \citet{HaEt10} highlight the relevance of the choice of the observation scale when estimating an occupancy distribution. In fact, as they clearly point out, the solution of MaxEnt depends on the scale on which the method is applied. They analyze a transformation $\lambda: A_{n_1,r} \longrightarrow A_{n_2,r}$ from a microscopic scale $n_1$ to a mesoscopic one $n_2$, where $n_1$ is a multiple of $n_2$, and in this respect a property of scale-consistency becomes relevant.\\
{%
In what follows we provide a formal definition of the concept of self-consistent distributions introduced by \citet{HaEt10}.
}%

\begin{definition}
Let $p^*_{n_1}$ and $p^*_{n_2}$ denote the solutions of MaxEnt at scales $n_1$ and $n_2$, respectively. We say that a probability distribution $p^*$ is \emph{scale-consistent} if
\[
p^*_{n_1} \circ \lambda^{-1} = p^*_{n_2} \,,
\]
where $p^*_{n_1} \circ \lambda^{-1}$ is the push-forward of  the measure $p^*_{n_1}$.
\end{definition}

Roughly speaking, a distribution is said to be scale consistent if solving MaxEnt at the scale $n_1$ and then applying the scale change to get an occupancy distribution on the scale $n_2$ coincides with solving MaxEnt directly at the coarse scale $n_2$.
The authors show that MB satisfies this property and BE does not. 
We can observe that the merging transformation $I_s$, as it has been defined in the previous section, does just coincide with the \emph{scale-transformation} in \cite{HaEt10}. Let us look, in fact, at the problem of detecting scale-consistent distributions from the perspective suggested by the arguments in Section~\ref{Sect:Merging}. In view of our results therein, we realize that their conclusions are a consequence of Proposition~\ref{Prop:composition:I} and Proposition~\ref{Prop:closure:I}. In this respect, we can conclude this brief discussion with the following simple result concerning the family $\mathcal{I}_{\mathcal{M}^{(a)}}$ defined in \eqref{def:class:IMa}.
 
\begin{proposition}
The probability distributions of the family $\mathcal{I}_{\mathcal{M}^{(a)}}$ are scale-consistent solutions of MaxEnt \eqref{MaxEnt:pb} over \mbox{$S = A_{n,r}$} under the constraints \eqref{MaxEnt:Ma:constraints}.
\end{proposition}

We feel that this analysis may deserve much more room and defer some more detailed work to future research.

\appendix

\section{Proofs}\label{Sect:Proofs}

\subsection{Proof of Proposition~\ref{Prop:composition:I}}
 
We want to show that the occupancy models $\left( I_{s_1} \circ I_{s_2} \right) (P)$ and $I_{s_1s_2} (P)$ are the same probability distribution in $\mathcal{P} \left( A_{n,r} \right)$; i.e. that the following diagram commutes:
\[
\xymatrix@R=3pt{
\mathcal{P} \left( A_{N,r} \right) \ar@{->}@//[r]^-{I_{s_2}} \ar@{->}@/_3pc/[rr]_{I_{s_1s_2}} & \mathcal{P} \left( A_{ns_1,r} \right) \ar@{->}@//[r]^-{I_{s_1}} & \mathcal{P} \left( A_{n, r} \right) 
}
\]
In the next table we fix some notation. \\
\begin{center}
\begin{tabular}{|cc|ccc|}
\hline
{\small \bf Model} & & \multicolumn{3}{c|}{\small \bf Occupancy Numbers} \\
\hline
&&&&\\
$P \in \mathcal{P} \left( A_{N,r} \right)$ & & $W_{j,h,k}$ & for & $\begin{array}{l} j=1, \dots, n \\ h=1, \dots, s_1 \\ k=1, \dots,s_2 \end{array}$ \\
&&&&\\
\hline
&&&&\\
$I_{s_2} (P) \in \mathcal{P} \left( A_{ns_1, r} \right)$ & & $Z_{j,h} = \displaystyle{\sum_{k=1}^{s_2} W_{j,h,k}}$ & for & $\begin{array}{l} j=1, \dots, n \\ h=1, \dots, s_1 \end{array}$ \\ 
&&&&\\
\hline
&&&&\\
$\left( I_{s_1} \circ I_{s_2} \right) (P) \in \mathcal{P} \left( A_{n,r} \right)$ & & $X_j = \displaystyle{\sum_{h=1}^{s_1} Z_{j,h}}$ & for & $j=1, \dots, n$ \\
&&&&\\
\hline
&&&&\\
$I_{s_1s_2} (P) \in \mathcal{P} \left( A_{n,r} \right)$ & & $X_j = \displaystyle{\sum_{h=1}^{s_1} \sum_{k=1}^{s_2} W_{j,h,k}}$ & for & $j=1, \dots, n$\\
&&&&\\
\hline
\end{tabular}
\end{center}
Now consider $P$ and apply transformations $I_{s_2}$ and $I_{s_1}$ sequentially. It yields
\begin{align}\label{composition:I:1}
P \left( X_1=x_1, \dots, X_n=x_n \right) &= \sum_{\mathbf{z} \in H_{n,r,\mathbf{x}}^{(s_1)}} P \left\{ Z_{1,1}=z_{1,1}, \dots, Z_{n,s_1}=z_{n,s_1} \right\} \nonumber\\
&= \sum_{\mathbf{z} \in H_{n,r,\mathbf{x}}^{(s_1)}} \; \sum_{\mathbf{w} \in H_{ns_1,r,\mathbf{z}}^{(s_2)}} P\left\{ W_{1,1,1}=w_{1,1,1}, \dots, W_{n,s_1,s_2}=w_{n,s_1,s_2} \right\} \,.
\end{align}
Alternatively, if we determine $I_{s_1s_2} (P)$, we get
\begin{equation}\label{composition:I:2}
P \left( X_1=x_1, \dots, X_n=x_n \right) = \sum_{\mathbf{w} \in H_{n,r,\mathbf{x}}^{(s_1s_2)}} P\left\{ W_{1,1,1}=w_{1,1,1}, \dots, W_{n,s_1,s_2}=w_{n,s_1,s_2} \right\} \,.
\end{equation}
The sets $H_{n,r,\mathbf{x}}^{(s_1)}$, $H_{ns_1,r,\mathbf{z}}^{(s_2)}$ and $H_{n,r,\mathbf{x}}^{(s_1s_2)}$ are defined similarly to \eqref{Constraints:set}. Observing that the sums \eqref{composition:I:1} and \eqref{composition:I:2} are comprised of the same terms concludes the proof.

\subsection{Proof of Proposition~\ref{Prop:closure:I}}
 
Let $Z_{1,1}, \dots, Z_{1,s}, Z_{2,1}, \dots, Z_{2,s}, \dots, Z_{n,1}, \dots, Z_{n,s}$ be the occupancy numbers corresponding to the $\mathcal{M}^{(a)}_{N,r}$-model $P$; meaning 
\[
P\left\{ Z_{1,1}=z_{1,1}, \dots, Z_{1,s}=z_{1,s}, \dots, Z_{n,1}=z_{n,1}, \dots, Z_{n,s}=z_{n,s} \right\} = \frac{\prod_{j=1}^n \prod_{h=1}^s a \left( z_{j,h} \right)}{C^{(a)}_{ns,r}} \,.
\]
We must prove that the vector $\left( X_1, \dots, X_n \right)$, where $X_j = \sum_{h=1}^{s} Z_{j,h}$ for $j=1, \dots, n$, is distributed as a $\mathcal{M}^{(a')}_{n,r}$-model for some suitable function $a'$. We have
\begin{align*}
P \{ X_1=x_1, \dots, X_n=x_n \} 
%
&= \sum_{\mathbf{z} \in H_{n,r,\mathbf{x}}^{(s)}} \frac{\prod_{j=1}^n \prod_{h=1}^s a \left( z_{j,h} \right)}{C^{(a)}_{ns,r}} \\
&= \frac{\prod_{j=1}^n \left[ \sum_{\mathbf{z}_j^{(s)} \in A_{s,x_j}} \prod_{h=1}^s a \left( z_{j,h} \right) \right]}{ C^{(a)}_{ns,r}} \\
&= \frac{\prod_{j=1}^n C^{(a)}_{s,x_j}}{C^{(a)}_{ns,r}} \,.
\end{align*}
Since $P \left\{ X_1=x_1, \dots, X_n=x_n \right\}$ is a probability distribution on $A_{n,r}$, it must be
\begin{equation}\label{id:norm:const}
C_{n,r}^{(a')} = C_{ns,r}^{(a)}
\end{equation}
with $a'(x)=C_{s,x}^{(a)}$ and the conclusion follows.

\subsection{Proof of Proposition~\ref{Prop:dropping:I}}

Let $X_1', \dots, X_n'$ denote the occupancy numbers of $K_1 (P') \in \mathcal{M}^{(a')}_{n,r-1}$. By Proposition~\ref{Prop:closure:I} we know $a' ( \cdot) = C_{s,\cdot}^{(a)}$ and moreover, since $a'$ satisfies \eqref{23}, Proposition~\ref{Prop:6.6} implies
\begin{equation}\label{drop:I:1}
P \left\{ X_1' = x_1', \dots, X_n'=x_n' \right\} = \frac{\prod_{j=1}^n C_{s,x_j'}^{(a)}}{C^{(a')}_{n,r-1}} \quad \mbox{ for } \mathbf{x}' \in A_{n,r-1} \,.
\end{equation}
Distribution \eqref{drop:I:1} can be obtained as merging of an occupancy model $\widetilde{P} \in \mathcal{M}^{(a)}_{ns,r-1}$, whose occupancy numbers $\widetilde{Z}_{1,1}, \dots, \widetilde{Z}_{n,s}$ are such that $\widetilde{\mathbf{Z}}^{(s)}_j \in A_{s,x_j'}$ for every $j=1, \dots, n$. This prove the first assertion in the statement. To conclude it remains to show that the following diagram commutes:
\[
\xymatrix{
\mathcal{M}^{(a)}_{ns,r}  \ar@{->}@//[d]_-{K_1} \ar@{->}@//[r]^-{I_s} & \mathcal{M}^{(a')}_{n,r} \ar@//[d]^-{K_1} \\
\mathcal{M}^{(a)}_{ns,r-1} \ar@{->}@//[r]_-{I_s} & \mathcal{M}^{(a')}_{n,r-1} 
}
\]
Since by hypothesis function $a$ fulfills \eqref{23}, by applying Proposition~\ref{Prop:6.6} we obtain that the distribution of $(Z'_{1,1}, \dots, Z'_{n,s})$, occupancy numbers of $K_1 (P) \in \mathcal{M}^{(a)}_{ns,r-1}$, is given by
\[
P \left\{ Z'_{1,1} = z'_{1,1}, \dots, Z'_{n,s} = z'_{n,s} \right\} = \frac{\prod_{j=1}^{n} \prod_{h=1}^s a(z'_{j,h})}{C^{(a)}_{ns,r-1}} \quad \mbox{ for } \mathbf{z}' \in A_{ns,r-1} \,.
\]
Now we merge groups of $s$ cells and we consider the occupancy model $\left( I_s \circ K_1 \right) (P) \in \mathcal{M}^{(a')}_{n,r-1}$; the joint distribution of its occupancy numbers $X_1, \dots, X_n$ is
\[
P \left\{ X_1 = x_1, \dots, X_n = x_n \right\} = \frac{\prod_{j=1}^n C_{s,x_j}^{(a)}}{C^{(a)}_{ns,r-1}} \quad \mbox{ for } \mathbf{x} \in A_{n,r-1} \,,
\]
which equals \eqref{drop:I:1} in view of the identity \eqref{id:norm:const}.

\subsection{Proof of Proposition~\ref{Prop:conditioning:I}}
 
We must show that $( K_{n,r}^{N,R} \circ I_s ) (P)$ and $( I_s \circ K_{ns,r}^{Ns,R} ) (P)$ are the same occupancy model over $A_{n,r}$; i.e. that the following diagram commutes: 
\[
\xymatrix{
\mathcal{M}^{(a)}_{Ns,R}  \ar@{->}@//[d]_-{K^{Ns,R}_{ns,r}} \ar@{->}@//[r]^-{I_s} & \mathcal{M}^{(a')}_{N,R} \ar@//[d]^-{K^{N,R}_{n,r}} \\
\mathcal{M}^{(a)}_{ns,r} \ar@{->}@//[r]_-{I_s} & \mathcal{M}^{(a')}_{n,r} 
}
\]
The proof relies on the closure property of the $\mathcal{M}^{(a)}$ class with respect to transformation $K^{N,R}_{n,r}$. Let $X_1, \dots, X_N$ denote the occupancy numbers of $P'=I_s(P)$ and $Z_{1,1}, \dots, Z_{N,s}$ those of $P \in \mathcal{M}^{(a)}_{Ns,R}$. \\
We start by determining the probability distribution corresponding to $K_{n,r}^{N,R} (P')$. Since $P'$ is a $\mathcal{M}^{(a')}_{N,R}$-model with $a'(\cdot)=C^{(a)}_{s,x}$, by Proposition~\ref{Prop:6.5} we get 
\begin{equation}\label{cond:I:via1}
P \left\{ X_1=x_1, \dots, X_n=x_n \left\vert \sum_{j=1}^n X_j = r \right. \right\} =  \frac{\prod_{j=1}^n C^{(a)}_{s,x_j}}{C^{(a')}_{n,r}} \quad \mbox{ for } \mathbf{x} \in A_{n,r}\,.
\end{equation}
On the other hand, if we consider $K_{ns,r}^{Ns,R} (P)$, by using Proposition~\ref{Prop:6.5} we obtain 
\begin{equation}\label{F1}
P \left\{Z_{1,1}=z_{1,1} \dots, Z_{n,s}=z_{n,s} \left\vert \sum_{j=1}^n \sum_{h=1}^s Z_{j,h} = r \right. \right\} = \frac{\prod_{j=1}^n \prod_{h=1}^s a(z_{j,h})}{C_{ns,r}^{(a)}} \quad \mbox{ for } \mathbf{z} \in A_{ns,r}
\end{equation}
and then, by applying $I_s$ it gives 
\begin{align}\label{cond:I:via2}
P \left\{ X''_1=x''_1, \dots, X''_n=x''_n \left\vert \sum_{j=1}^n X''_j = r \right. \right\} 
%
&\stackrel{\mbox{\tiny \eqref{F1}}}{=} \sum_{\mathbf{z} \in H_{n,r,\mathbf{x}''}^{(s)}} \frac{\prod_{j=1}^n \prod_{h=1}^s a \left( z_{j,h}\right)}{C^{(a)}_{ns,r}} \nonumber\\
&\hspace{.1cm} = \frac{\prod_{j=1}^n C^{(a)}_{s,x''_j}}{C^{(a)}_{ns,r}} \quad \mbox{ for } \mathbf{x}'' \in A_{n,r}\,,
\end{align}
where $X''_1, \dots, X''_n$ are the occupancy numbers of a $\mathcal{M}^{(a')}_{n,r}$-model. The conclusion follows by observing that \eqref{cond:I:via1} and \eqref{cond:I:via2} are equal, since \eqref{id:norm:const} holds.

\subsection{Proof of Proposition~\ref{Prop:MaxEnt:sol}}

The proof relies on the technique of Lagrange multipliers. Let $\lambda_0, \lambda_1, \dots, \lambda_m$ be $m+1$ multipliers. The Lagrangian associated with problem \eqref{MaxEnt:pb} is given by
\[
\mathcal{L} (p,\lambda_0, \lambda_1, \dots,\lambda_m) = - \sum_{s \in S} p(s) \log p(s) - \lambda_0 \left( \sum_{s \in S} p(s) - 1 \right) - \sum_{k=1}^m \lambda_k \left( \sum_{s \in S} f_k (s) p(s) - c_k \right). 
\]
Solutions of \eqref{MaxEnt:pb} correspond to the critical points of $\mathcal{L}$. By solving $\nabla \mathcal{L} = 0$, it yields
\[
p^*(s) = \frac{e^{- \sum_{k=1}^m \lambda_k f_k(s)}}{\sum_{s \in S} e^{- \sum_{k=1}^m \lambda_k f_k(s)}} \,,
\]
where $\lambda_1, \dots, \lambda_m$ must be chosen so that $p^*$ fulfills the constraints $\sum_{s \in S} f_k (s) p(s) = c_k$, for all $k=1, \dots, m$. If the constraints cannot be satisfied for any values of $\lambda$'s, then the maximum entropy distribution does not exist.\\
Uniqueness of the MaxEnt solution follows from the convexity of entropy together with the linearity of the constraints.

\subsection{Proof of Proposition~\ref{Prop:MaxEnt:Ma}}

By adapting the general solution \eqref{Gibbs:distr:gen} to our specific case, we obtain
\[
p^* (\mathbf{x}) = \frac{\prod_{j=1}^n a^{-\lambda} (x_j)}{\sum_{\mathbf{x} \in A_{n,r}}\prod_{j=1}^n a^{-\lambda} (x_j)} \quad \mbox{ for } \mathbf{x} \in A_{n,r}\,.
\]
Therefore, if we take constant $c$ so that the Lagrange multiplier $\lambda=-1$, the distribution solving MaxEnt is \eqref{def:M(a)-model} as wanted.

\section*{Acknowledgements}
The authors are very grateful to Dr. Marco Formentin for suggesting reference \citet{HaEt10}, for fruitful conversations and insights on the maximum entropy principle. 

\end{document}